\documentclass[a4,11pt]{article}
\usepackage[a4paper,top=43mm, bottom=18mm, left=25mm, right=25mm]{geometry} 
\usepackage{amsmath,amssymb}
\usepackage{amsthm, amsfonts}
\usepackage{mathrsfs}
\usepackage{url}
\usepackage{fancyhdr}
\usepackage{paralist}
\usepackage{graphics} 
\usepackage{epsfig} 
\usepackage{graphicx} 
\usepackage{epstopdf}
\usepackage[colorlinks=true]{hyperref}
\usepackage{subcaption}
\usepackage[english]{babel}
\usepackage{color} 
\usepackage{pstricks}
\usepackage{verbatim}
\usepackage{multirow}
\usepackage{amstext} 
\usepackage{bm}

\usepackage{floatrow}  
\usepackage{wrapfig}

 \usepackage{setspace}
\hypersetup{urlcolor=blue, citecolor=red}

\newtheorem{Remark}{Remark}
\newtheorem{Definition}{Definition}
\newtheorem{Proposition}{Proposition}

\newcommand{\R}{\mathbb{R}}
\newcommand{\diag}{\text{diag}}
\newcommand{\1}{\mathbb{1}}

\def\ds{\displaystyle}

\DeclareSymbolFont{msbm}{U}{msb}{m}{n}

\newcommand{\bfx}{{\bf x}}

\newcommand{\bfu}{{\bf u}}

\newcommand{\bfF}{{\bf F}}

\newcommand{\msL}{{\mathcal L}}



\begin{document}





\title{A Stability notion for the viscous Shallow Water Discrete-Velocity Boltzmann Equations
 }

\author{ Mapundi K. Banda\footnote{Department of Mathematics and Applied Mathematics, Botany Building 2 - 10, \newline University of Pretoria, Hatfield 0028, South Africa.\newline Tel:+27 12 420 2544 \hspace{25pt} fax:+27 12 420 3893 \newline
\emph{E-mail address}: \texttt{mapundi.banda@up.ac.za} (Mapundi Banda)}\\
Tumelo R.A.~Uoane\footnote{CSIR Modelling and Digital Sciences, Meiring Naude Road, Brummeria,
Pretoria, 0001, South Africa. \newline \emph{E-mail address}: \texttt{tumelo.uoane@gmail.com} (Tumelo Uoane)}
}

\date{}                                          
\maketitle

\begin{abstract}
The theoretical stability of Lattice Boltzmann Equations modelling Shallow Water Equations in the special case of reduced gravity is investigated. A stability notion as applied in incompressible Navier-Stokes equations in Banda, M. K., Yong, W.- A. and Klar, A: A stability notion for lattice Boltzmann equations. SIAM J. Sci. Comput. {\bf 27(6)}, 2098-2111 (2006) is used. It is found that to maintain stability a careful choice of the value of the reduced gravity must be made. The stability notion is employed to investigate different shallow water lattice Boltzmann models. Results are tested using the Lattice Boltzmann Method for various values of the governing parameters of the flow. It is observed that even for the discrete model the reduced gravity has a significant effect on the stability.
\end{abstract}


\vspace{.5cm}



\section{Introduction}\label{section1}

In this paper, a stability notion for the viscous Shallow Water Lattice Boltzmann  equations (SWLBE) is discussed. The Shallow Water Equations are popular in modelling flow phenomena which includes flood waves, dam breaks, tidal flows in an estuary and coastal water regions, and bore wave propagation in rivers. Herewith a few of these will be highlighted: wind-driven ocean circulation \cite{Salmon,ZhongFengGao2005}, three-dimensional planetary geostrophic equations \cite{Salmon R}, or the atmospheric circulation of the northern hemisphere with ideal boundary conditions \cite{FengZhaoTsutaharaJi2002}. Such equations are derived from the depth-averaged incompressible Navier-Stokes equations and usually they include continuity and momentum equations.
For such real-life processes, it is imperative that the numerical approaches applied to simulate the flow are based on accurate and efficient models. Above all the stability of such models must also be classified. 

In general the lattice Boltzmann Method (LBM) is based on a special discretization of Boltzmann-type kinetic equations, a system of hyperbolic equations with stiff source terms \cite{Sterling}. The hyperbolic equations are a statistical physics formulation of fluid flow. It is an approach based on the description of the flow of distribution functions (at the mesoscopic scale) of fluid particles in discrete space instead of the classical description based on macroscopic variables. This formulation is applied on fluid flow which in the macroscopic scale models shallow water flow \cite{Salmon,Dellar,Guido,Zhou} which is the limit of slow varying solutions. The basic idea is to replace the nonlinear differential equations of macroscopic fluid dynamics by a simplified description modeled on the kinetic theory of gases. The advantage of this kinetic-type approach is that the advection terms are linear but the local source terms are stiff. The linearity can be exploited to simplify programming or simulations for complex geometry, irregular topography, structured meshes while the stiff source terms are treated using local operators. It has also been known to be effective for implementation on parallel computer architectures \cite{kandhai1998}.

In general the stability of the continuous kinetic models is known. These satisfy a dissipative entropy condition (Boltzmann's H-theorem) \cite{Cercignani}. The same can not be said about the discrete-velocity models, the reader may refer to the discussion in \cite{Yong A,Yong C}. Instead the Lattice Boltzmann equations have been constructed to satisfy some physical requirements like Galilean invariance and isotropy, to possess a velocity-independent pressure and no compressible effects \cite{Koelman JMVA,Zou}. Furthermore,  alternative stability conditions have been developed. These include: the structural stability in \cite{Yong B}, the sub-characteristic condition \cite{Jin} and the dissipative entropy principles \cite{Bouchut}.

In the following work, stable LB models will be identified by using stability conditions in \cite{Yong B,Banda}. In most of the models used, it is not yet rigorously proven that the diffusive limit of the discrete-velocity Boltzmann equation are SWEs at least in the regime of smooth flow. But we can remark that incompressible fluids are modelled using either SWEs or the N-S equations. The latter satisfies the diffusive limit of the discrete-velocity Boltzmann equation, see \cite{Junk}, when certain models are used. These models are similar to the ones used in this work. Therefore, it is reasonable to consider the stability condition as a new requirement in constructing LB equations for the SWEs. A previous discussion on the stability of the shallow water lattice Boltzmann models was presented in \cite{Dellar}. There in a von Neumann approach was applied. To further the discussion, in this paper an alternative notion \cite{Banda} will be used to demonstrate the stability structure of the discrete-velocity models.

It should be pointed out that the stability theory used here is different from the previous works \cite{Dellar} on the stability of the lattice Boltzmann method (an explicit difference scheme). In \cite{Dellar} stability analysis was based on the von Neumann stability analysis and the resulting growth matrix was not treated analytically.  In contrast, the theory presented here is based on a rigorous asymptotic analysis \cite{Yong B} for the lattice Boltzmann equations (partial differential equations). This analysis was applied to the lattice Boltzmann equations for incompressible Navier-Stokes in \cite{Banda}. A linearised stability of the lattice Boltzmann method (the completely discrete form) was presented in \cite{Junk2}. There in a few examples of lattice Boltzmann methods for which the structural hypothesis holds were presented. Thus it is useful not only for the lattice Boltzmann method but also for other discretizations of the hyperbolic systems. Moreover, the derivations of the parameter relations is purely analytic (see Section \ref{Determination of Parameters}). It must be emphasized that only two-dimensional models are considered.

The popularly used reduced gravity model is discussed in Section \ref{section2} which also briefly discusses an existence result. To explain how the stability requirement guides the construction of the LB equations, we will show that the LB models for the SWEs are stable using Definition (\ref{def1}) in Section \ref{section3}. We will do so by testing the stability structure on some examples which will be shown in Section \ref{section4}. In other models, we will also investigate the parameter range for which the models are stable. Computational experiments were undertaken on examples which are used commonly in literature, to confirm the applicability of the stability structure and the results are presented in Section \ref{section4}.

\section{The Shallow Water Models and the Discrete-Velocity Formulation}\label{section2}

\subsection{The Shallow Water Model}\label{ShallowWater}

The two-dimensional shallow water equations including friction and Coriolis forces take the form: \begin{eqnarray}
\partial_t{h} + \partial_x(hu_1) + \partial_y(hu_2) &=& 0,\nonumber \\
\partial_t{(hu_1)} + \partial_x\left(hu_1^2+\frac{1}{2}gh^2\right) + \partial_y\left(hu_1u_2\right) &=& - gh\partial_x{Z} + \nabla\cdot(h\nu\nabla{u_1}) + \nonumber \\
&& \frac{1}{\rho_0}\bigl(\mathcal{T}_{wx} - \mathcal{T}_{bx}\bigr) - \Gamma hu_2,\label{swe2}\\
\partial_t{(hu_2)} + \partial_x\left(hu_1u_2\right) + \partial_y\left(hu_2^2+\frac{1}{2}gh^2\right) &=& - gh\partial_y{Z} + \nabla\cdot(h \nu\nabla{u_2}) +\nonumber \\
&& \frac{1}{\rho_0}\bigl(\mathcal{T}_{wy} - \mathcal{T}_{by}\bigr) + \Gamma hu_1,\nonumber
\end{eqnarray}
At the macroscopic level, the water depth, $h$, and depth-averaged water velocity $\bfu = (u_1,u_2)^T$ are obtained from solving the shallow water equations in Equation \eqref{swe2}. In this equation $u_1(x,y,t)$ and $u_2(x,y,t)$ are the depth-averaged water velocity in $x$- and $y$-direction, $\rho_0$ is the water density, $g$ is the gravitational acceleration, $Z$ is the bottom topography, $\nu$ is the horizontal kinematic viscosity, $\Gamma$ is the Coriolis parameter defined by $\Gamma = 2\omega\sin\phi$ (where $\omega = 0.000073\;\text{rad}\;\text{s}^{-1}$ is the angular velocity of the earth and $\phi$ the geographic latitude), and $\nabla = (\partial_x,\partial_y)^T$ is the gradient operator. The bottom stresses $\mathcal{T}_{bx}$ and $\mathcal{T}_{by}$ are the bed shear stresses in the $x$- and $y$-direction, respectively, defined with respect to the depth-averaged velocities as
\begin{equation}
\mathcal{T}_{b_x} = \rho_0C_bu_1\sqrt{u_1^2 + u_2^2},\qquad \mathcal{T}_{b_y} = \rho_0C_bu_2\sqrt{u_1^2 + u_2^2},\label{slopef}
\end{equation}
where $C_b$ is the bed friction coefficient, which may be either constant or estimated as $C_b = {g}/{C_z^2}$. Note that $C_z = h^{1/6}/n_b$ is the Chezy constant, in which $n_b$ is the Manning roughness coefficient at the bed. The surface stresses $\mathcal{T}_{w_x}$ and $\mathcal{T}_{w_y}$ are wind stresses defined using the wind velocity,
\begin{equation}
\mathcal{T}_{w_x} = \rho_0C_w w_1\sqrt{w_1^2 + w_2^2},\qquad \mathcal{T}_{w_y} = \rho_0C_w w_2\sqrt{w_1^2 + w_2^2},\label{wind}
\end{equation}
where $C_w$ is the coefficient of wind friction and ${\bf w} =(w_1,w_2)^T$ is the velocity of the wind at $10\;m$ above the water surface. It is usually defined by \cite{Bermudez A}
\[C_w = \rho_a\left(0.75+0.067\sqrt{w_1^2 + w_2^2}\right)\times10^{-3},\]
where $\rho_a$ is the air density. Note that other coefficients of wind friction in \eqref{wind} can also be applied.

It has to be pointed out that it is well known that the shallow water problems \eqref{swe2} can be derived from the depth-averaged incompressible Navier-Stokes equations with the assumption that the vertical scale is much smaller than any typical horizontal scale and the pressure is hydrostatic. Thus the quantity $gh$ defines the geopotential. 

\begin{Remark}\label{existence} In \cite{Sundbye} an existence proof for the Dirichlet problem for viscous shallow water equations excluding $Z$, the bottom topography, and the Coriolis forces was given. The existence proof gives guidance  to our choice of the equilibrium values for stability analysis in Section \ref{section3}. In summary under certain assumptions, it was proved that the viscous shallow water equation has a unique global solution in time and a unique equilibrium state $(\bar{h}, {\mathbf 0})$ with initial conditions \begin{equation} \label{Sund:init} \bfu(x, y, 0) = \bfu_0(x, y), \qquad h(x, y, 0) = h_0(x, y), \quad \text{ for } (x, y) \in \Omega \end{equation}
and Dirichlet boundary conditions \begin{equation} \label{Sund:bc} \bfu(x, y, t) = {\mathbf 0}, \qquad \text{ for } (x, y) \in \partial\Omega, \; t \ge 0. \end{equation}


For a detailed theorem and proof the reader may refer to \cite{Sundbye}. Further in \cite{Sundbye} the requirement for the positivity of the fluid height for $t \ge 0$ was also established. 
\end{Remark}

In the next section, the lattice Boltzmann equations for the shallow water flow equations \eqref{swe2} are presented. A discussion of the discrete-velocity model will also be briefly discussed.  
\subsection{Discrete-Velocity Boltzmann equation for Shallow Water Flows}
\label{Discrete-velocity Boltzmann equation for Shallow Water Flow with source terms}
The continuum two-dimensional kinetic Equation \eqref{kinetic} is considered
\begin{equation}
\frac{\partial f}{\partial t} + \bm{\xi} \cdot \nabla f = J(f) + F.\label{kinetic}
\end{equation}
Equation \eqref{kinetic} describes the evolution of a particle density $f(\bfx,\bm{\xi},t)$ with $\bfx = (x,y) \in \R^2$ the spatial variable  and $\bm{\xi} = (\xi_1,\xi_2) \in \R^2$ are the microscopic velocities of the particle distribution $f$. In \eqref{kinetic}, $J$ is the collision term, and $F$ is the effect of external forces. The left hand side of Equation \eqref{kinetic} represents the linear transport of fluid particles.

\begin{figure}
\centering\includegraphics[width=0.35\textwidth,height=0.25\textheight,angle=0]{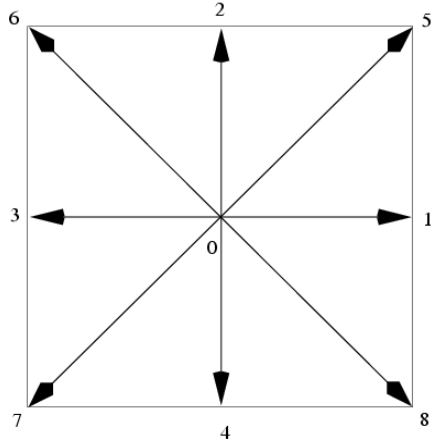}
\caption{Links in the D2Q9 lattice Boltzmann method.}\label{figure2}
\end{figure}

For the discrete-velocity models in two space dimensions, assume
\[\bm{\xi} \in \{\bm{\xi}_0, \bm{\xi}_1, \ldots, \bm{\xi}_{N-1}\},\]
with $\bm{\xi}_i \in \R^2$. Here, the D2Q9 square lattice model \cite{Qian} as sketched in Figure \ref{figure2} is an example of the discrete-velocity model, with the velocity vectors of particles defined by
\begin{equation*}
\renewcommand{\arraystretch}{1.1}
\begin{aligned}
\bm{\xi}_0&=\begin{pmatrix} 0\\0 \end{pmatrix},& \bm{\xi}_1&=\begin{pmatrix} 1\\0 \end{pmatrix},& \bm{\xi}_2&=\begin{pmatrix} 0\\1\end{pmatrix}, & \bm{\xi}_3&=\begin{pmatrix}-1\\0 \end{pmatrix}, & \bm{\xi}_4&=\begin{pmatrix} 0\\-1 \end{pmatrix},\\[2ex]
 \bm{\xi}_5&=\begin{pmatrix} 1\\1 \end{pmatrix},& \bm{\xi}_6&=\begin{pmatrix} -1\\1 
\end{pmatrix}, & \bm{\xi}_7&=\begin{pmatrix} -1\\-1 \end{pmatrix}, & \bm{\xi}_8&=\begin{pmatrix} 1\\-1 \end{pmatrix}.
\end{aligned}
\end{equation*}

In the discrete-velocity case, the $\bm{\xi}$-dependence of the particle distribution $f(\bfx, \bm{\xi}, t)$ is determined through $N$ functions 
\[f_i(\bfx,t) = f(\bfx, \bm{\xi}_i, t),\qquad i = 0, 1, \ldots, N-1.\]
Hence the discrete-velocity equation can be written as:
\begin{equation}
\label{eqn3.1}
\frac{\partial f_{i}}{\partial t}+ \bm{\xi}_{i}.\nabla f_{i} = J_{i}(f)\quad(i = 0,1,\ldots,N-1).
\end{equation}

The physical variables, the water depth, $h$, and the velocity $\bfu$, are defined in terms of the distribution function as
\begin{equation}
h(\bfx,t) = \sum_i f_i(\bfx, t), \qquad h\bfu = \sum_i \bm{\xi}_if_i(\bfx,t).\label{equib}
\end{equation}
In most approaches for the lattice Boltzmann applications, the collision operator $J(f)$ in \eqref{kinetic} is of BGK-type \cite{bgk}
\begin{equation}
J(f) = - \frac1{\tau}(f - f^{eq}), \label{eqn:bgk}
\end{equation}
where the parameter $\tau > 0$ is called the relaxation time and $f^{eq}$ is the equilibrium distribution. In the shallow water case, $f^{eq}$ depends on $f$ through the parameters $h$ and $\bfu$ which are calculated according to \eqref{equib}. The local equilibrium function satisfies the following conditions
\begin{equation}
\sum_i f_i^{eq} = h,\qquad \sum_i \bm{\xi}_i f_i^{eq} = h\bfu,\qquad \sum_i \bm{\xi}_i \bm{\xi}_i f_i^{eq} = P(h){\bf I} + h \bfu \otimes \bfu,\label{feq3}
\end{equation}
where $P(h) = \frac12 gh^2$ such that the lattice Boltzmann equation approaches the solution of the two-dimensional shallow water equations. In \eqref{feq3}, $\mathbf{I}$ denotes the $2\times2$ identity matrix. 
For the standard D2Q9-model with nine velocities, $f^{eq}$ takes the form  \cite{Salmon, Dellar}
\begin{equation}
f_i^{eq}(h, \bfu) = \begin{cases} 
\ds h - f_0^*h\left(\frac{15}{2}gh - \frac{3}{2}\bfu^2\right), & i = 0, \\[3.5ex]
\ds f_i^*h\left(\frac{3}{2}gh + 3\bm{\xi}_i\cdot \bfu + \frac{9}{2}(\bm{\xi}_i\cdot \bfu)^2 - \frac{3}{2} \bfu^2\right), & i = 1,\dots,8,
\end{cases}\label{feqi}
\end{equation}
with the D2Q9 weight factors 
\begin{equation}
f_i^* = \begin{cases} 
\ds \frac49, & i = 0, \\[1.5ex]
\ds \frac19, & i = 1,2,3,4, \\[1.5ex]
\ds \frac1{36}, & i = 5,6,7,8.
\end{cases}\label{feqi2}
\end{equation}

To obtain the macroscopic equations from equation \eqref{kinetic}, the Chapman-Enskog asymptotic expansion can be employed \cite{Salmon,Dellar}. The LB equation \eqref{kinetic} with equilibrium function \eqref{feqi} and collision term \eqref{eqn:bgk} results in the solution of the SWE \eqref{swe2} with a force term ${\mathbf{F}}$:
\begin{equation}
\renewcommand{\arraystretch}{2}
\bfF(\bfx,t) = \begin{pmatrix}\ds - gh\partial_x Z + \frac 1 {\rho_0}(\mathcal{T}_{w_x}-\mathcal{T}_{b_x})- \Gamma hu_2 \\
- gh\partial_y Z + \frac 1 {\rho_0}(\mathcal{T}_{w_y}-\mathcal{T}_{b_y})+ \Gamma hu_1 
\end{pmatrix},
\end{equation}as required. Thus, the external force terms such as wind stress, Coriolis force, and bottom friction are easily included in the model by introducing them into the force term ${\mathbf{F}}$. For details on this multi-scale expansion, the reader may refer to  \cite{Salmon,Dellar,ZhongFengGao2005}.

Hence, using a special discretization of the above BGK approximation \cite{Salmon,Zhou}, the following fully discrete lattice Boltzmann equation is obtained
\begin{equation}
f_i(\bfx + \bm{\xi}_i\Delta x, t + \Delta t) - f_i(\bfx, t) = -\frac{\Delta t}{\tau}\left(f_i - f^{eq}_i\right) +3\Delta t f_i^* \bm{\xi}_{i}\cdot\bfF(x,t),\label{feq4}
\end{equation}
where $\Delta t$ is the time scale, $\Delta x$ is the reference length. A stability analysis for such a discrete form for incompressible Navier-Stokes Equations was presented in \cite{Junk2}. A similar analysis for the shallow water equation is not yet available. 

By applying a Taylor expansion on equation \eqref{feq4} and the Chapman-Enskog procedure, it can be shown that the solution of the discrete lattice Boltzmann equation \eqref{feq4} with the equilibrium function \eqref{feqi} results in the solution of the shallow water equations \eqref{swe2} with a lattice Boltzmann viscosity, $\hat\nu$, defined as
\begin{equation}
{\hat\nu} = \frac{1}{6}\Bigl(2\hat\tau - 1\Bigr)\label{lbmvisco}
\end{equation}
where $\hat\tau = \tau/\Delta{t}$ is the scaled relaxation time. This viscosity is related to the physical viscosity in \eqref{swe2} 
by the relation
\begin{equation}
\frac{\nu}{{\hat\nu}} = e^2\Delta{t},
\end{equation}
where $e = \Delta x/\Delta t$ denotes the velocity along a unit link \cite{Salmon,Guido,Zhou}. Further, for the limit of small Mach number which is of interest here $e^2 \ll gh$ for consistency \cite{Salmon,Dellar,ZhongFengGao2005}.

\section{Stability Structure}\label{Stability Structure}\label{section3}

In this section the LB equation \eqref{eqn3.1} derived from a particular discretization to obtain a $d$-dimensional, $N$-velocity Boltzmann equation is considered:
\begin{Definition}{Stability Structure\cite{Junk2, Banda}}:
\\Let $f_{*}$ be a constant state satisfying $J(f_{*}) = 0$. The Lattice Boltzmann Equation \eqref{eqn3.1} is called stable at $f = f_{*}$ if there is an invertible matrix $P \in \R^{N\times N}$ such that $P^TP$ is diagonal, $\diag(a_1, a_2, \ldots, a_N)$, and $$PJ_f(f_{*}) = - diag(\lambda_1, \lambda_2, \ldots, \lambda_N)P$$ with $\lambda_i = 0$ for $i \leq d + 1$ and $\lambda_i > 0$ for $i > d + 1$. Here $J_f(f)\in \R^{N\times N}$ is the Jacobian of $J(f) = (J_{1}(f), J_{2}(f), \ldots, J_{N}(f))^T$. In this case the lattice Boltzmann Equation is said to be stable at $f = f_{*}$. The triple $(P, a, \lambda)$ is referred to as the stability structure at $f = f_{*}$.

\label{def1}
\end{Definition}

In the above definition, $d$ represents the space dimension, $J_f(f_{*})$ is the Jacobian matrix and $\lambda_i$ is an eigenvalue.
\begin{Remark}\label{remark3.1}
\begin{itemize}
\item[(a)] This definition is based on the stability conditions \cite{Junk2,Yong B, Banda} for hyperbolic systems with source terms. 
\item[(b)] In general, a consistent and stable lattice Boltzmann model implies convergence. It is hoped that this can be proven using the approach in \cite{Banda}. In the case of shallow water equations, consistency is not yet proven and it is not within the scope of this work.\label{remark3.2}
\end{itemize}
\end{Remark} 

The stability structure introduced above will be verified for the example models below. The models to be used are taken from \cite{Salmon, Zhou}.

\subsection{The Stability Structure for the D2Q7 Model}
Next we consider D2Q7-velocity model, in which
\[\bm{\xi}_0  = (0,0),\]
\[\{\bm{\xi}_i: i \in \{1,\ldots, 6\}\} = \left\{e\left[\mbox{cos}\left(\frac{(i - 1)\pi}{3}\right), \mbox{sin}\left(\frac{(i - 1)\pi}{3}\right)\right]\right\}.\]
 The collision terms  are
\[J_{i}(f) = \frac{f^{eq}_{i}(h,\bm{u}) - f_{i}}{\tau}\]
where
\begin{equation}
\label{eqn3.12}
f^{eq}_{i}(h,\bm{u}) = \left\{\begin{array}{clcr}
h - \displaystyle{\frac{g h^2}{e^2} + \frac{h\textbf{u}^2}{e^2}},& \mbox{$i=0$} \\ \\
\displaystyle{\frac{g h^2}{6e^2} + \frac{h \bm{\xi}_i \textbf{u}}{3 e^2} +  \frac{2 h (\bm{\xi}_i\cdot\textbf{u})^2}{3 e^4} - \frac{h \textbf{u}^2}{2 e^2}},& \mbox{$i \in \{1,\ldots, 6\}$}.
\end{array}\right.\
\end{equation}  with
\[h = \sum_{i = 0}^{6}f_i, \qquad  h\textbf{u} = \sum_{i = 0}^{6}\bm{\xi}_i f_i.\]

Let us assume the following (these are straight-forward to verify directly)
\begin{equation}
\label{eqn3.13}
\begin{array}{clcr}
\displaystyle{\sum_{i = 0}^{N}f^{eq}_i = h = \sum_{i = 0}^{N}f_i},\qquad
\displaystyle{\sum_{i = 0}^{N}\bm{\xi}_i f^{eq}_i = h \textbf{u} = \sum_{i = 0}^{N}\bm{\xi}_i f_i.}
\end{array}
\end{equation}
Deriving the Jacobian from Equation (\ref{eqn3.12}), gives
\begin{equation}
\label{eqn3.14}
\frac{\partial f^{eq}_i (h ,\textbf{u})}{\partial f_j} = \left\{\begin{array}{clcr}
1 - \displaystyle{\frac{2gh}{e^2} + \frac{2 \bm{\xi}_j \textbf{u}}{e^2} - \frac{\textbf{u}^2}{e^2}},& \mbox{$i = 0$;}\\ \\
\displaystyle{\frac{gh}{3e^2} + \frac{\bm{\xi}_j \bm{\xi}_i}{3e^2} +  \frac{4\bm{\xi}_j \bm{\xi}_i\cdot\textbf{u}}{3e^4} - \frac{2(\bm{\xi}_i\cdot\textbf{u})^2}{3e^4} - \frac{\bm{\xi}_j\textbf{u}}{e^2} + \frac{\textbf{u}^2}{2e^2}},& \mbox{$i \neq 0$}.
\end{array}\right.\
\end{equation}
By using Equation (\ref{eqn3.13}), we deduce from Equation (\ref{eqn3.14}) that
\[[f^{eq}_{f_{i}}(h,\textbf{u})]^2 = [f^{eq}_{f_{i}}(h,\textbf{u})],\]
that is, the Jacobian $[f^{eq}_{f_{i}}(h,\textbf{u})]$ is a projection matrix. Thus, the eigenvalues of 
\[J_f(h ,\textbf{u}) = ([f^{eq}_{f_{i}}(h,\textbf{u})] - I_7) / \tau\]
are 0 and $\displaystyle{-\frac{1}{\tau}}$. Take $f_{*} = f^{eq}(\bar{h},\bm{0})$ then
\begin{equation}
\label{eqn3.15}
\frac{\partial f^{eq}_i (\bar{h},\bm{0})}{\partial f_j} = \left\{\begin{array}{clcr}
\displaystyle{\frac{e^2 - 2 g\bar{h}}{e^2}},& \mbox{$i =  0$} \\ \\
\displaystyle{\frac{g\bar{h}}{3e^2} + \frac{\bm{\xi}_j \bm{\xi}_i}{3e^2}},& \mbox{$i \neq 0$}.
\end{array}\right.\
\end{equation}
Let $\bm{c} = (1,1,\ldots,1) \in \R^7$, $\bm{\xi} = (\bm{\xi}_0,\bm{\xi}_1, \ldots, \bm{\xi}_6)$. Also let 
\begin{equation}
\label{eqn3.17}
B_0 = \mbox{diag} \left[\frac{e^2}{e^2 - 2 g\bar{h}}, \frac{3e^2}{g\bar{h}}\textbf{I}_6\right].
\end{equation} 
such that 
\begin{equation}\label{eqn:symmetry}
B_0[f^{eq}_{f_{i}}(h,\textbf{u})] = \frac1{3e^2} \bm{\xi}^T\bm{\xi}  + B_0\Bigl(\frac{e^2 - 2g\bar{h}}{e^2}, \frac{g\bar{h}}{3e^2},\ldots, \frac{g\bar{h}}{3e^2}\Bigr)\bm{c},
\end{equation}

From the above choice of $B_0$, we need to choose $g$ such that, $B_0$ remains positive definite. Therefore, we set 
\[g < \frac{e^2}{2\bar{h}} \hspace{0.5cm} \mbox{and} \hspace{0.5cm} g \neq 0.\]

We deduce from Equation \eqref{eqn:symmetry} that the rank of $[f^{eq}_{f_{i}}(\bar{h},\bm{0})]$ is 3. Since $[f^{eq}_{f_{i}}(h,\textbf{u})]$ is a projection matrix and 
\[\tau J_f(f_{*}) = [f^{eq}_{f_{i}}(\bar{h},\bm{0})] - I_7\]
then, the rank of $J_f(f_{*})$ is 4. 

On the other hand, since $B_0 [f^{eq}_{f_{i}}(\bar{h},\bm{0})]$ is symmetric and $B_0$ is symmetric  positive definite, it is well known that there is an invertible matrix $P$ such that 
\[ B_0 = P^T P \hspace{1cm} \mbox{and} \hspace{1cm} B_0 \tau J_f(f_{*}) = P^T \Lambda P\]
with $\Lambda$ a diagonal matrix. We may as well assume that 
\[\Lambda = -\mbox{diag}(0, 0, 0, 1, 1, 1, 1).\]
Thus we have proven,
\begin{Proposition}
If $g < \frac{e^2}{2\bar{h}}$ then the 2-dimensional 7-velocity model is stable at $f_{*} = f^{eq}(\bar{h},\bm{0})$.
\label{prop2}
\end{Proposition}

\begin{Remark}
The lattice Boltzmann model used above was developed by Zhou \cite{Zhou} using the $7$-speed hexagonal lattice. The model was developed in the same manner to that of the $9$-speed square lattice. 
\label{remark3.5}
\end{Remark}

\subsection{The Stability Structure for the D2Q9 Model}\label{Determination of Parameters}
In the following section, some parameters will be fixed for D2Q9 LB models. By doing so, it can be assumed that the LB models are stable for those fixed parameters. The models to be used are taken from \cite{Salmon,Dellar}. The following examples are used:
 Consider D2Q9-velocity model, with 
\[\bm{\xi}_0 = (0,0),\]
\[\{\bm{\xi}_i: i = 1,2,3,4\} = \{(\pm{e},0)^T, (0,\pm{e})^T\},\]
\[\{\bm{\xi}_i: i = 5,6,7,8\} = \{(\pm{e},\pm{e})^T\}.\]

Further the following moments are listed:
\[h = \sum_{i = 0}^{8}f_i,\quad h \textbf{u} = \sum_{i = 0}^{8}\bm{\xi}_i f_i, \quad \Pi = \sum_{i = 0}^{8}\bm{\xi}_i \bm{\xi}_i f_i.\]
To show its stability, firstly, the following are assumed (it is straightforward to verify these)
\begin{eqnarray}
\label{eqn3.24}
&& \sum_{i = 0}^{8}f^{eq}_i = h = \sum_{i = 0}^{8}f_i,
\sum_{i = 0}^{8}\bm{\xi}_i f^{eq}_i = h \textbf{u} = \sum_{i = 0}^{8}\bm{\xi}_i f_i, \nonumber \\
&& \\
 \nonumber
&& \Pi^{(eq)} = \sum_{i = 0}^{8}\bm{\xi}_i \bm{\xi}_i f^{eq}_i
 = P(h)\textbf{I} + h \textbf{uu}.
\end{eqnarray}
see \cite{Dellar}. The equilibrium distribution, the so called Salmon's equilibrium \cite{Salmon}, is given as:
\begin{equation}
\label{eqn3.26}
f^{eq}_i (h,\textbf{u}) = \left\{ \begin{array}{clcr}
 h - \displaystyle{\frac{5gh^2}{6e^2} - \frac{2}{3e^2}h \textbf{u}^2} & \mbox{ $i = 0$}\\ \\
\displaystyle{\frac{gh^2}{6e^2} + \frac{1}{3e^2} h \bm{\xi}_i\cdot\textbf{u} + \frac{1}{2e^4}h(\bm{\xi}_i\cdot\textbf{u})^2 - \frac{1}{6e^2}h \textbf{u}^2},& \mbox{$1\leq i \leq 4$}\\ \\
\displaystyle{\frac{gh^2}{24e^2} + \frac{1}{12e^2} h \bm{\xi}_i\cdot\textbf{u} + \frac{1}{8e^4}h(\bm{\xi}_i\cdot\textbf{u})^2 - \frac{1}{24e^2}h \textbf{u}^2},& \mbox{$5\leq i \leq 8$.}
\end{array}\right.\
\end{equation}

Deriving
\begin{equation}
\label{eqn3.27}
\frac{\partial f^{eq}_i (h,\textbf{u})}{\partial f_j} = \begin{cases}
1 - \displaystyle{\frac{5}{3e^2}gh - \frac{4}{3e^2} \bm{\xi}_j\cdot\textbf{u} + \frac{2}{3e^2}\textbf{u}^2}, & \mbox{ $i = 0$}\\
& \\
\displaystyle{\frac{1}{3e^2}gh + \frac{1}{3e^2}\bm{\xi}_j\bm{\xi}_i + \frac1{e^4}(\bm{\xi}_i\cdot\bm{\xi}_j)(\bm{\xi}_i\cdot\bm{u}) } & \\
\ds{- \frac{1}{2e^4}(\bm{\xi}_i\cdot\bm{u})^2 - \frac{1}{3e^2}\bm{\xi}_j\cdot\textbf{u} + \frac{1}{6e^2}\textbf{u}^2},& \mbox{$1\leq i \leq 4$}\\
& \\
\displaystyle{\frac{1}{12e^2}gh + \frac{1}{12e^2}\bm{\xi}_j \bm{\xi}_i + \frac{1}{4e^2}(\bm{\xi}_i\cdot\bm{\xi}_j)(\bm{\xi}_i\cdot\textbf{u})} & \\
\ds{ - \frac{1}{8e^4}(\bm{\xi}_i \cdot\textbf{u})^2 - \frac{1}{12e^2}\bm{\xi}_j\cdot\textbf{u} + \frac{1}{24e^2}\textbf{u}^2},& \mbox{$5\leq i \leq 8$.}
\end{cases}
\end{equation}
By using Equation (\ref{eqn3.24}), we deduce from Equation (\ref{eqn3.27}) that
\[[f^{eq}_{f_{i}}(h,\textbf{u})]^2 = [f^{eq}_{f_{i}}(h,\textbf{u})],\]
that is, the Jacobian $[f^{eq}_{f_{i}}(h,\textbf{u})]$ is a projection matrix. Thus, the eigenvalues of 
\[J_f(h,\textbf{u}) = ([f^{eq}_{f_{i}}(h,\textbf{u})] - I_9) / \tau\]
are 0 and $\displaystyle{-\frac{1}{\tau}}$. Take $f_{*} = f^{eq}(\bar{h},\bm{0})$, then
\begin{equation}
\label{eqn3.28}
\frac{\partial f^{eq}_i (\bar{h},\bm{0})}{\partial f_j} = \left\{ \begin{array}{clcr}
 1 - \displaystyle{\frac{5}{3e^2}g\bar{h}} & \mbox{$i = 0$}\\ \\
\displaystyle{\frac1{3e^2}g\bar{h} + \frac{1}{3e^2}\bm{\xi}_j\cdot \bm{\xi}_i} ,& \mbox{$1\leq i \leq 4$}\\ \\
\displaystyle{\frac{1}{12e^2}g\bar{h} + \frac{1}{12e^2}\bm{\xi}_j\cdot \bm{\xi}_i },& \mbox{$5\leq i \leq 8$.}
\end{array}\right.\
\end{equation}
Let $\bm{c} = (1,1,\ldots,1) \in \R^9$, $\bm{\xi} = (\bm{\xi}_0,\bm{\xi}_1, \ldots, \bm{\xi}_8)$. Also let 
\begin{equation}\label{eqn3.31}
C_0 = \frac{3e^2}{g\bar{h}} \hspace{0.07cm}\mbox{diag} \left[\frac{g\bar{h}}{3e^2 - 5g\bar{h}}, \textbf{I}_4, 4\textbf{I}_4\right].
\end{equation}
Then we get
\begin{align}
\label{eqn3.29}
C_0 \left[\frac{\partial f^{eq}_i (\bar{h},\bm{0})}{\partial f_j}\right] = \frac1{g\bar{h}}\bm{\xi}^T\bm{\xi} + \nonumber \\ C_0\Bigl(\frac{3e^2-5g\bar{h}}{3e^2}, \frac{g\bar{h}}{3e^2}, \frac{g\bar{h}}{3e^2}, \frac{g\bar{h}}{3e^2}, \frac{g\bar{h}}{3e^2}, \frac{g\bar{h}}{12e^2}, \frac{g\bar{h}}{12e^2}, \frac{g\bar{h}}{12e^2}, \frac{g\bar{h}}{12e^2}\Bigr)^T\bm{c},
\end{align}
The square matrix $C_0$ needs to be symmetric and positive definite. For this to hold using Equation (\ref{eqn3.28}), we see that
\[1 - \displaystyle{\frac{5}{3e^2}} g\bar{h} > 0 \hspace{0.5cm} \mbox{and} \hspace{0.5cm} g \neq 0\]
i.e, $g \in (0, \displaystyle{\frac{3e^2}{5\bar{h}}})$.
Hence the right hand side of Equation \eqref{eqn3.29} is symmetric. The rank of $[f^{eq}_{f_{i}}(\bar{h},\textbf{u})] |_{\bm{u} = \bm{0}}$ is 3. Since $[f^{eq}_{f_{i}}(h,\textbf{u})]$ is a projection matrix and 
\[\tau J_f(f_{*}) = [f^{eq}_{f_{i}}(\bar{h},\textbf{u})] |_{\bm{u} = \bm{0}} - I_9\]
then, the rank of $J_f(f_{*})$ is 6. 

On the other hand, since $C_0 [f^{eq}_{f_{i}}(\bar{h},\textbf{u})] |_{\bm{u} = \bar{0}}$ is symmetric and $C_0$ is symmetric  positive definite, it can be concluded that there is an invertible matrix $P$ such that 
\[ C_0 = P^T P \hspace{1cm} \mbox{and} \hspace{1cm} C_0 \tau J_f(f_{*}) = P^T \Lambda P\]
with $\Lambda$ a diagonal matrix. Therefore, we may take $\Lambda$ to be:
\[\Lambda = -\mbox{diag}(0, 0, 0, 1, 1, 1, 1, 1, 1).\]
Hence we have proved,

\begin{Proposition}
The 2-dimensional 9-velocity model with \eqref{eqn3.26} is stable at $f_{*} = f^{eq}(\bar{h},\bm{0})$ if $g \in (0, \displaystyle{\frac{3e^2}{5\bar{h}}})$.
\label{prop3}
\end{Proposition}


\subsection{The Stability Structure for the D2Q9 Model with Parameter $\lambda$}
Consider another D2Q9-velocity model which was investigated in \cite{Dellar}, with

\begin{equation}
\label{eqn3.35}
f^{eq}_i (h,\textbf{u}) = \begin{cases}
\displaystyle{\frac{(8 + \lambda)}{9}h - \frac{(4 + \lambda)}{6e^2}gh^2 - \frac{2}{3e^2}h \textbf{u}^2},  & \mbox{ $i = 0$}\\ 
& \\
\displaystyle{\frac{(1 - \lambda)}{18}h + \frac{(1 + \lambda)}{12e^2}gh^2 + \frac{1}{3e^2} h \bm{\xi}_i\cdot \textbf{u}} & \\
\ds{ + \frac{1}{2e^4}h(\bm{\xi}_i\cdot\textbf{u})^2 - \frac{1}{6e^2}h \textbf{u}^2},& \mbox{$1\leq i \leq 4$}\\ 
& \\
\displaystyle{\frac{(\lambda - 1)}{36}h + \frac{(2 - \lambda)}{24e^2}gh^2 + \frac{1}{12e^2} h \bm{\xi}_i\cdot\textbf{u}} & \\
\ds{ + \frac{1}{8e^4}h(\bm{\xi}_i\cdot\textbf{u})^2 - \frac{1}{24e^2}h \textbf{u}^2},& \mbox{$5\leq i \leq 8$.}
\end{cases}
\end{equation}

The Jacobian of the above equilibrium distribution takes the form:


\begin{equation}
\label{eqn3.36}
\frac{\partial f^{eq}_i (h,\textbf{u})}{\partial f_j} = \begin{cases}
 \displaystyle{\frac{(8 + \lambda)}{9} - \frac{(4 + \lambda)}{3e^2}gh - \frac{4}{3e^2} \bm{\xi}_j\cdot\textbf{u}} & \\
 \ds{ + \frac{2}{3e^2}\textbf{u}^2}, & \mbox{ $i = 0$}\\ 
 & \\
 \displaystyle{\frac{(1 - \lambda)}{18} + \frac{(1 + \lambda)}{6e^2}gh + \frac{1}{3e^2}\bm{\xi}_j\cdot\bm{\xi}_i} & \\
 \ds{ + \frac1{e^4}(\bm{\xi}_i\cdot\bm{\xi}_j)(\bm{\xi}_i\cdot\bm{u})} & \\
 \ds{ - \frac{1}{2e^4}(\bm{\xi}_i\cdot\bm{u})^2 - \frac{1}{3e^2}\bm{\xi}_j\cdot\textbf{u} + \frac{1}{6e^2}\textbf{u}^2},& \mbox{$1\leq i \leq 4$}\\
 & \\
 \displaystyle{\frac{(\lambda - 1)}{36} + \frac{(2 - \lambda)}{12e^2}gh + \frac{1}{12e^2}\bm{\xi}_j \bm{\xi}_i} & \\
  \ds{ + \frac{1}{4e^2}(\bm{\xi}_i\cdot\bm{\xi}_j)(\bm{\xi}_i\cdot\textbf{u})}  & \\
 \ds{- \frac{1}{8e^4}(\bm{\xi}_i \cdot\textbf{u})^2  - \frac{1}{12e^2}\bm{\xi}_j\cdot\textbf{u} + \frac{1}{24e^2}\textbf{u}^2},& \mbox{$5\leq i \leq 8$.}
\end{cases}
\end{equation}

By using (\ref{eqn3.24}), we deduce from (\ref{eqn3.36}) that
\[[f^{eq}_{f_{i}}(h,\textbf{u})]^2 = [f^{eq}_{f_{i}}(h,\textbf{u})]\]
that is, the Jacobian $[f^{eq}_{f_{i}}(h,\textbf{u})]$ is a projection matrix. Thus, the eigenvalues of $J_f(h,\textbf{u}) = ([f^{eq}_{f_{i}}(h,\textbf{u})] - I_9) / \tau$ are 0 and $\displaystyle{-\frac{1}{\tau}}$. 
Take $f_{*} = f^{eq}(\bar{h},\bm{0})$, then
\begin{equation}
\label{eqn3.37}
\frac{\partial f^{eq}_i (\bar{h},\bm{0})}{\partial f_j} = \left\{ \begin{array}{clcr}
\displaystyle{\frac{(8 + \lambda)}{9} - \frac{(4 + \lambda)}{3e^2}g\bar{h}} ,& \mbox{ $i = 0$}\\ \\
\displaystyle{\frac{(1 - \lambda)}{18} + \frac{(1 + \lambda)}{6e^2}g\bar{h} + \frac{1}{3e^2}\bm{\xi}_j\cdot\bm{\xi}_i},& \mbox{$1\leq i \leq 4$}\\ \\
\displaystyle{\frac{(\lambda - 1)}{36} + \frac{(2 - \lambda)}{12e^2}g\bar{h} + \frac{1}{12e^2}\bm{\xi}_j\cdot \bm{\xi}_i},& \mbox{$5\leq i \leq 8$.}
\end{array}\right.\
\end{equation}
Let $\bm{c} = (1,1,\ldots,1) \in \R^9$, $\bm{\xi} = (\bm{\xi}_0,\bm{\xi}_1, \ldots, \bm{\xi}_8)$, $\bm{c}_4 = (1,1,\ldots,1) \in \R^4$. Also let 
\begin{multline}
\label{eqn3.39}
D_0 = \mbox{diag} \left[\Big(\frac{9e^2}{(8 + \lambda)e^2 - 3(4 + \lambda) g\bar{h}}\Big), \Big(\frac{18e^2}{(1 - \lambda)e^2 + 3(1 + \lambda){g}\bar{h}}\Big)\textbf{I}_4,\right. \\
\left.\Big(\frac{36e^2}{(\lambda - 1)e^2 + 3(2 - \lambda)g\bar{h}}\Big)\textbf{I}_4\right].
\end{multline} 

Using Equation \eqref{eqn3.39} one obtains
\begin{equation}
\label{eqn3.38}
D_0 \left[\frac{\partial f^{eq}_i (\bar{h},\bm{0})}{\partial f_j}\right] = \Gamma \bm{\xi}^T\bm{\xi}  + D_0 \Psi^T\bm{c}
\end{equation}
where 
\begin{eqnarray*}
\Psi &=& \left(\left(\frac{(8 + \lambda)e^2 - 3(4 + \lambda) g\bar{h}}{9e^2}\right), \left(\frac{(1 - \lambda)e^2 + 3(1 + \lambda){g}\bar{h}}{18e^2}\right)\bm{c}_4,\right. \\
&&\left. \left(\frac{(\lambda - 1)e^2 + 3(2 - \lambda)g\bar{h}}{36e^2}\right)\bm{c}_4\right);\\
\Gamma &=& \mbox{diag} \left[\left(\frac{9e^2}{(8 + \lambda)e^2 - 3(4 + \lambda) g\bar{h}}\right), \left(\frac{6e^2}{(1 - \lambda)e^2 + 3(1 + \lambda){g}\bar{h}}\right)\textbf{I}_4,\right. \\
&& \left. \left(\frac{3e^2}{(\lambda - 1)e^2 + 3(2 - \lambda)g\bar{h}}\right)\textbf{I}_4\right];
\end{eqnarray*}
in which the right hand side of Equation \eqref{eqn3.38} is a symmetric matrix if
\begin{equation}
\label{eqn3.41}
\frac{6e^2}{(1 - \lambda)e^2 + 3g\bar{h} + 3 \lambda{g}\bar{h}} = \frac{3e^2}{(\lambda - 1)e^2 + 6 g\bar{h} - 3\lambda g\bar{h}},
\end{equation} which is true if $\lambda = 1$.

We need to choose the parameters $g$ and $\lambda$ such that $D_0$ is positive definite. In particular, we notice that the model in (\ref{eqn3.35}) is similar to the model in (\ref{eqn3.26}) for the above value of $\lambda$. Then, for the stability structure [\ref{def1}] to hold, $g \in (0, \displaystyle{\frac{3e^2}{5\bar{h}}})$.

From Equation \eqref{eqn3.38} it can be deduced that the rank of $[f^{eq}_{f_{i}}(\bar{h},\textbf{u})] |_{\bm{u} = \bm{0}}$ is 3. Since $[f^{eq}_{f_{i}}(h,\textbf{u})]$ is a projection matrix and 
\[\tau J_f(f_{*}) = [f^{eq}_{f_{i}}(\bar{h},\textbf{u})] |_{\bm{u} = \bm{0}} - I_9\]
then, the rank of $J_f(f_{*})$ is 6. 

On the other hand, since $D_0 [f^{eq}_{f_{i}}(\bar{h},\textbf{u})] |_{\bm{u} = \bm{0}}$ is symmetric and $D_0$ is symmetric  positive definite, it can be concluded that there is an invertible matrix $P$ such that 
\[ D_0 = P^T P \hspace{1cm} \mbox{and} \hspace{1cm} D_0 \tau J_f(f_{*}) = P^T \Lambda P\]
with $\Lambda$ a diagonal matrix. We may as well assume that 
\[\Lambda = -\mbox{diag}(0, 0, 0, 1, 1, 1, 1, 1, 1).\]
We also observe that when $g = \displaystyle{\frac{e^2}{3\bar{h}}}$ in Equation (\ref{eqn3.37}) the parameter $\lambda$ is arbitrary. 

Thus it can be proved that,

\begin{Proposition}
If $\lambda = 1$ and $\ds{g \in (0, \frac{3e^2}{5\bar{h}})}$ or if $\lambda$ is arbitrary and $\ds{g = \frac{e^2}{3\bar{h}}}$, then the D2Q9-velocity model with \eqref{eqn3.35} is stable at $f_{*} = f^{eq}(\bar{h}, \bm{0})$.
\label{prop4}
\end{Proposition}

\begin{Remark}\label{remark3.9}
\begin{itemize}
\item[(a)] In \cite{Dellar}, it was claimed that the parameter $\lambda$ was adjustable to give a positive equilibrium in the state at rest, i.e, when $\textbf{u} = 0$. The stability structure (\ref{def1}) above gave only two relationships, when $\lambda = 1$ and arbitrary, proving the claim in \cite{Dellar}.
\item[(b)] To complete the discussion on stability structure, we discuss a D2Q5 model that was presented in \cite{Salmon}. Unlike other LB models, this model has no momentum advection term. The stability structure of this model was also investigated. It was found that the model is stable at $f_{*} = f^{eq}(\bar{h},\bm{0})$ if  $\ds{g < \frac{e^2}{2\bar{h}}}$.
\end{itemize}
\end{Remark}

From the above examples, it was shown that the stability requirement can be regarded as a reasonable guide for a good choice of parameters. When the choice of parameters do not satisfy the stability condition (\ref{def1}), unstable results might be obtained. Numerical tests will be presented in the next section to verify these results.

\section{Numerical Results}\label{section4}

In this section, the stability criteria as discussed in Section \ref{Determination of Parameters} by Propositions \ref{prop3} and \ref{prop4} is tested numerically. Three examples that are widely applied in literature are considered. These are: the steady flow over a hump \cite{Goutal}, tidal wave flow \cite{Bermudez A}, and flow over a sudden-expansion channel \cite{ZhouJG}. The main goal is to show that when reasonable ranges for the corresponding parameters are chosen, stable results are obtained. Alternatively, when choices of parameters are outside the suggested range in the propositions, then stable results can not be guaranteed. The accuracy of the Lattice Boltzmann  Method (stability structure [\ref{def1}]), is also demonstrated by comparing the numerical predictions with analytical solutions.

\subsection{Example 1: Steady flow over a hump}
In this example, the convergence in time towards the steady flow over the hump is shown. This example was also considered by the working group on dam break modelling \cite{Goutal} and used in \cite{Vazquez-Cendon ME} to test an upwind discretization for the bed slope source term. 

\begin{figure}[ht]
\centering
\includegraphics[width=0.55\textwidth,height=0.3\textheight,angle=0]{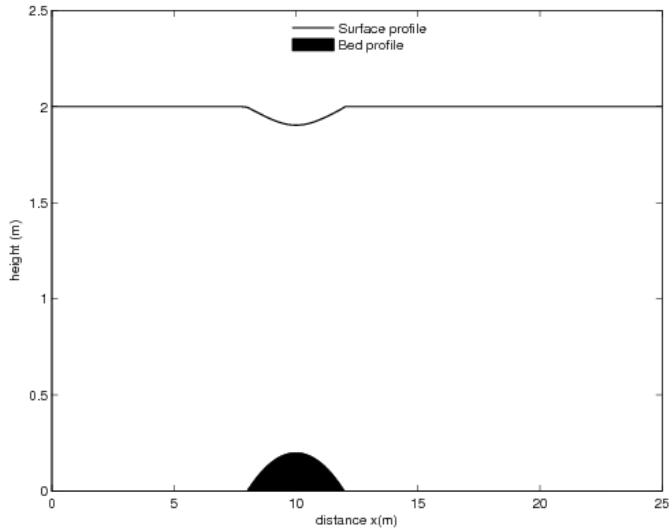}
	\caption{Steady subcritical flow over a hump: Illustration of the profile of water surface and bottom.}
	\label{hump1}
\end{figure}

A one-dimensional steady flow in a $25$m long and $1$m wide channel with a hump is defined by
\begin{equation*}
z_{b}(x) = \left\{ \begin{array}{clcr}
 &0.2 - 0.05(x - 10)^2, \hspace{1cm}\mbox{if}& \mbox{$8\leq x \leq 12$};\\ \\
&0, & \mbox{otherwise.}
\end{array}\right.\
\end{equation*}
The initial conditions are given by
\[h(x,0) = 2 \hspace{0.1cm}\mbox{m} - z_b(x) \hspace{0.5cm}\mbox{and}\hspace{0.5cm} u(x,0) = 0\hspace{0.1cm} \mbox{m/s}\]
as illustrated in Figure \ref{hump1}. We expect to observe that for steady subcritical flow passing over the hump on a bed slope, the water surface over the hump drops. The analytical solution is given in \cite{Goutal}. 

This example is used as a test problem to verify that Propositions \ref{prop3} and \ref{prop4} hold, starting with the former.
The following channel boundary conditions were prescribed, the water level $h = 2$ m is used as the outflow boundary condition and the discharge $q = 4.42$ $m^2/s$ is imposed at the inflow boundary; the slip or non-slip boundary conditions are used at the solid walls. In the lattice Boltzmann implementation, for the no-slip condition, the bounce-back scheme is used and for slip conditions, a zero gradient of the distribution function normal to the solid wall is used. The lattice speed $e = 15$ m/s and $\tau = 1.5$ are also used. 

We define the global relative error $R$ by 
\begin{equation}
\label{eqn5.9}
R = \sqrt{\sum_{i}{\left(\frac{h_{i}^{n} - h_{i}^{n-1}}{h_{i}^{n}} \right)^2}},
\end{equation}
as defined in \cite{Zhou}. The $h_{i}^{n}$ and $h_{i}^{n-1}$ represent the local water depth at the current and previous time levels, respectively. For the scheme to converge to a steady solution, the convergence criterion is taken as $R < 5 \times 10^{-6}$.

The three lattice sizes, $125\times50$, $250\times50$ and $500\times50$ which correspond to $\Delta x = 0.2$ m, $\Delta x = 0.1$ m and $\Delta x = 0.05$ m are used in the initial computations to test their effects on lattice solutions. For numerical computation the gravitational acceleration  $g$ ranges from 0.006 and 0.09, i.e $g\in(0.006, 0.09)$. The choice of $g$ for numerical computation was motivated by the fact that $g \in (0,\frac{3}{5})$ and computed on the lattice speed $e = 15$ m/s, giving $g \in (0,\frac{3}{5e})$. Steady state solutions were obtained from different values of $g$ used in the computation, refer to Table \ref{num:gravity}. 
\begin{table}[t]
	\caption{The summary of gravity for different lattices, "$-$" implies that there was no convergence.}
	\centering
	\label{num:gravity}
		\begin{tabular}{c c c}
			\hline\noalign{\smallskip} 
			\textbf{Gravity (g)} & \textbf{Lattice sizes} & \textbf{Number of iterations}\\ [3pt]
			\hline \hline
			\noalign{\smallskip}
	  	\multirow{3}{*}{0.09}	&	$125\times50$  & |\\								
																&	$250\times50$  & |\\
																&	$500\times50$  & |\\
			\hline													
			\multirow{3}{*}{0.07}	&	$125\times50$  & 19513\\							
																&	$250\times50$  & | \\
																&	$500\times50$  & | \\
			\hline 
		 \multirow{3}{*}{0.03}	  &	$125\times50$  & 19873\\
																&	$250\times50$  & 39170\\
																&	$500\times50$  & |\\
			\hline
			\multirow{3}{*}{0.009}	&	$125\times50$  & 21333\\
																&	$250\times50$  & 40034\\
																&	$500\times50$  & 59700\\
			\hline
			\multirow{3}{*}{0.006}	&	$125\times50$  & 24165\\
																&	$250\times50$  & 40319\\
																&	$500\times50$  & 60048\\  
			\noalign{\smallskip}													
			\hline\hline
		\end{tabular}
\end{table}
When values of $g$ outside the required range were used, the method was unstable. For example, when $g = 0.07$ steady state solution is reached only at $125 \times 50$ lattice points and the solution does not converge when the grid is refined. On the other hand when gravity ($g$) is reduced, better results are obtained (when $g = 0.03$ and 0.006). 

\begin{figure}[ht]
\centering\includegraphics[width=0.55\textwidth,height=0.3\textheight,angle=0]{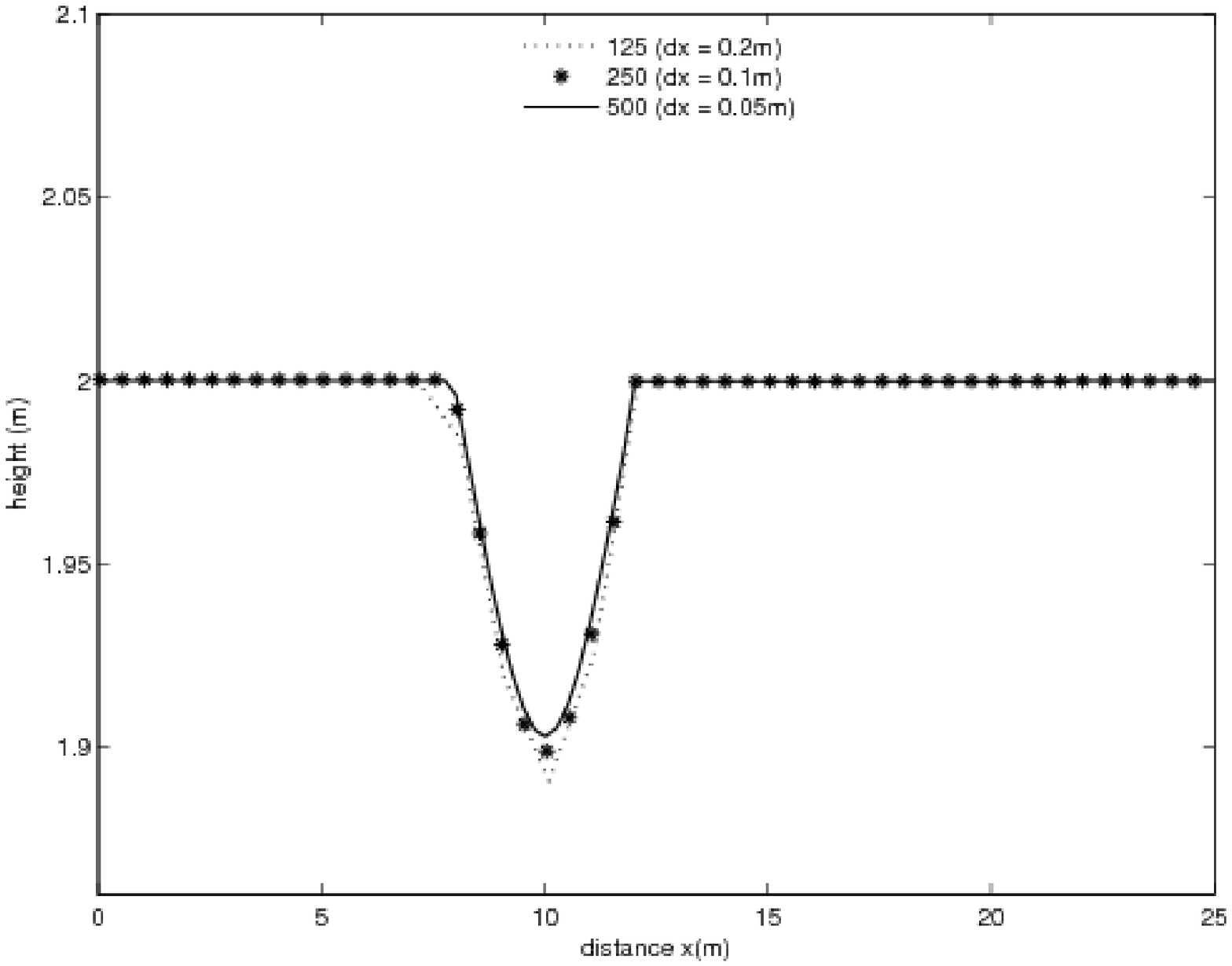}
	\caption{Steady subcritical flow over a hump: Effect of the lattice size on solutions.}
	\label{Three_lattices}
\end{figure}

The value of $g = 0.009$ was chosen when comparing numerical results between different lattice sizes. There was little difference found, refer to Figure \ref{Three_lattices}. 

The results further indicate that, for the values satisfying the stability structure, when lattice sizes become smaller better results are obtained, i.e. the results of $\Delta x = 0.1$ m and $\Delta x = 0.05$ m are almost the same, but there is a small difference between $\Delta x = 0.2$ m and $\Delta x = 0.1$ m. Hence the results at $\Delta x = 0.05$ m are preferred.

\begin{figure}[ht]
\centering\includegraphics[width=0.55\textwidth,height=0.3\textheight,angle=0]{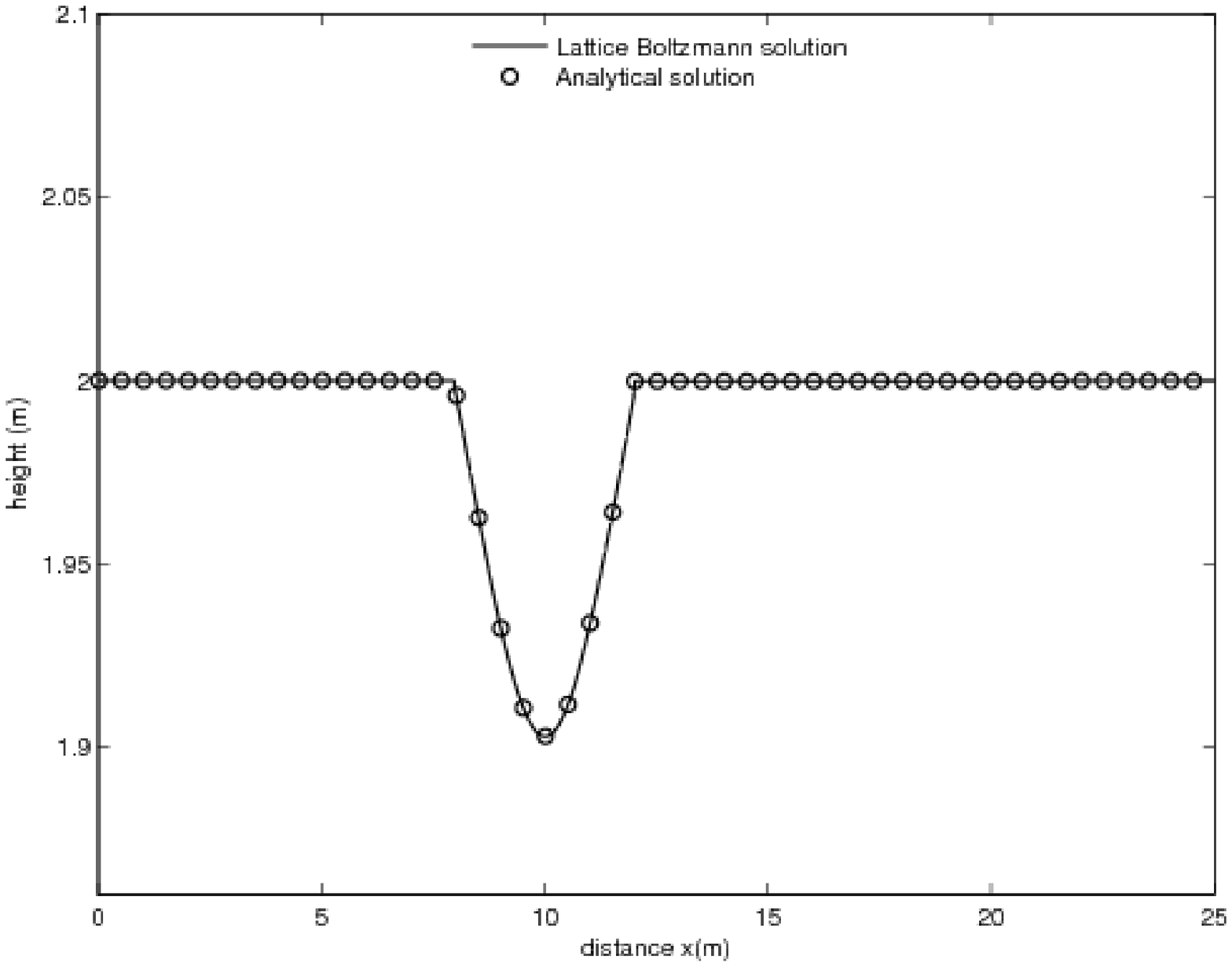}
	\caption{Steady subcritical flow over a hump: Comparison of the water surface.}
	\label{hump_surface}
\end{figure}

The accuracy of the approach was tested by comparing the computed steady water surface with the analytical solution as depicted in Figure \ref{hump_surface}, showing an excellent agreement.

\begin{figure}[ht]
\centering\includegraphics[width=0.55\textwidth,height=0.3\textheight,angle=0]{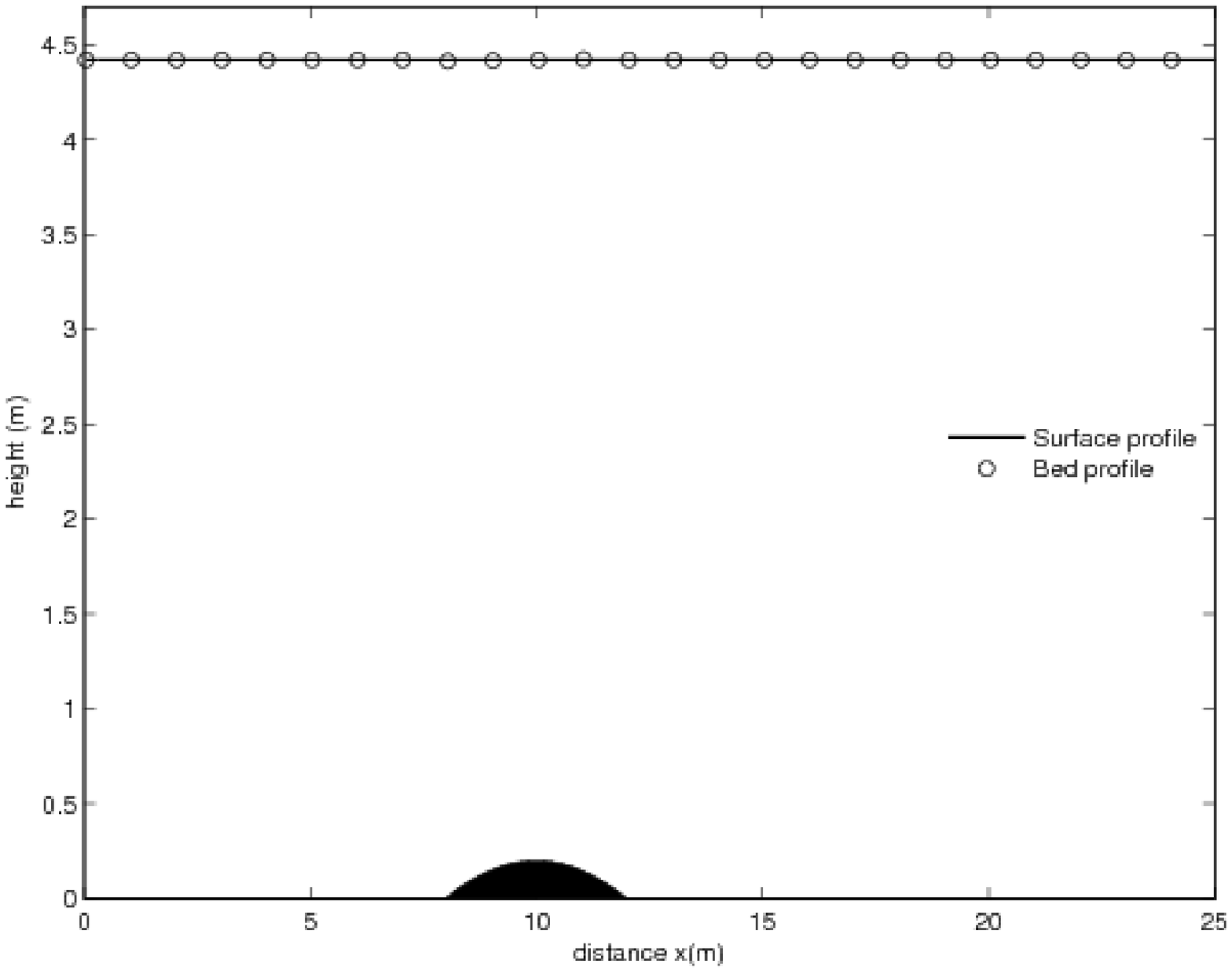}
	\caption{Steady subcritical flow over a hump: Comparison of discharge.}
	\label{discharge}
\end{figure}

 The $L^2$- error norm was used to verify the results, defined as
\begin{equation}
\label{eqn:error}
\left\|\textbf{c}\right\|_{L^2} = \sqrt{\frac{\sum_{ij}{|c^n - \tilde{c}(x_i,y_j,t_n)|^2}}{\sum_{ij}{|\tilde{c}(x_i,y_j,t_n)|^2}}},
\end{equation}
where $c^n$ is the computed LB solution and $\tilde{c}(x_i,y_j,t_n)$ is the analytical solution, respectively, at time $t_n$ and lattice point $(x_i,y_j)$. It was found that, the comparison of the computed LB solution with the analytical solution indicates that the relative error for the water depth is 0.325 $\%$. To test the conservative property of the model, the numerical solution of the discharge was computed and is depicted in Figure \ref{discharge}. The relative error was about 0.18 $\%$.

 This suggests that the model is conservative and accurate. Note that, the above results were based on $\Delta x = 0.05$ m lattice size.

To check if Proposition \ref{prop4} holds, the parameter $\lambda$ was varied from $-6$ to $14$ using $g = \frac{1}{3e}$ on different lattice sizes, refer to Table \ref{num:lambda}. It is interesting to note that when the value of $\lambda$ increases in magnitude, leads to unstable results. It was also shown in \cite{Dellar} that when distribution functions change signs, it leads to a stable equilibrium distribution function for the SWEs. 
\begin{table}[t]
	\caption{The summary of the value of $\lambda$ for different lattices using $g = \frac{1}{3e}$.}
	\centering
	\label{num:lambda}
		\begin{tabular}{c c c}
			\hline\noalign{\smallskip} 
			\textbf{Parameter $\lambda$} & \textbf{Lattice sizes} & \textbf{Number of iterations}\\ [3pt] 
			\hline \hline
			\noalign{\smallskip}
	  	\multirow{3}{*}{-6}	  &	$125\times50$  & 21333\\								
																&	$250\times50$  & 40034\\
																&	$500\times50$  & 59700\\
			\hline													
			\multirow{3}{*}{0}   	&	$125\times50$  & 21333\\								
																&	$250\times50$  & 40034\\
																&	$500\times50$  & 59700\\
			\hline 
		 \multirow{3}{*}{3}	    &	$125\times50$  & 21333\\
																&	$250\times50$  & 40034\\
																&	$500\times50$  & 59700\\
			\hline
			\multirow{3}{*}{6.7}	  &	$125\times50$  & 21333\\
																&	$250\times50$  & 40034\\
																&	$500\times50$  & |\\												
			\hline\hline
			\noalign{\smallskip}\hline
		\end{tabular}
\end{table}
To conclude the stability structure has been demonstrated that it can be used in the choice of parameters in order to obtain stable numerical simulations. It is believed that the instability of very large $\lambda$'s is rather a numerical artefact. Hence there is need to undertake numerical analysis on the full discrete lattice Boltzmann method itself to further investigate this artefact.

\subsection{Example 2: Tidal wave flow\cite{Bermudez A} }

\begin{figure}[ht]
\centering
\includegraphics[width=0.55\textwidth,height=0.3\textheight,angle=0]{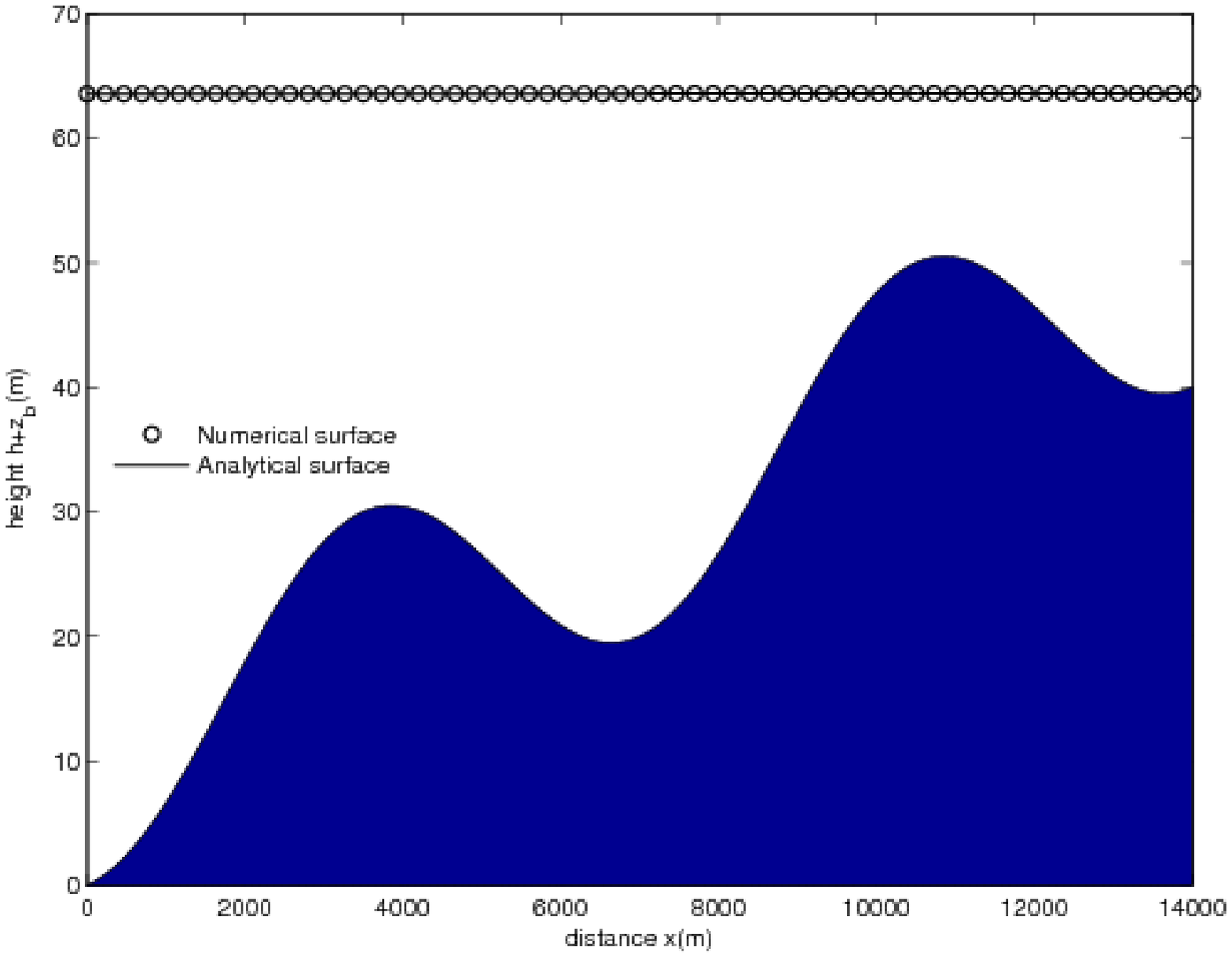}
	\caption{Numerical and analytical free surface for the tidal wave flow at time $t = 9117.5$ s.}
	\label{level}
\end{figure}

In this example a one-dimensional problem of a tidal wave in a channel was considered. In \cite{Bermudez A} this problem was used to test an upwind discretization of the bed slope source. The following is the description of the problem: the bed topography is defined by (refer to Figure \ref{level})
\begin{equation*}
H(x) = 50.5 - \frac{40 x}{L} + 10 \mbox{sin} \left( \pi \left( \frac{4 x}{L} - \frac{1}{2} \right) \right),
\end{equation*}
where $L = 14$ km is the length of the channel and $H(x)$ is the partial depth between a fixed reference level and the bed surface, giving $z_b = H(0) - H(x)$. 

The initial conditions for the water height and velocity are given by
\begin{equation*}
h(x,0) = H(x), \qquad u(x,0) = 0.
\end{equation*}

At the inflow and outflow of the channel, the water height and velocity, respectively, is defined by
\begin{equation*}
h(0,t) = H(0) + 4 - 4 \mbox{sin} \left( \pi \left( \frac{4t}{86400} - \frac{1}{2} \right) \right), \qquad u(L,t) = 0.
\end{equation*}
In  \cite{Bermudez A}, the asymptotic analytical solution for this test example was given by
\begin{equation}
\label{eqn5.7}
h(x,t) = H(x) + 4 - 4 \mbox{sin} \left( \pi \left( \frac{4t}{86400} - \frac{1}{2} \right) \right)
\end{equation}
and
\begin{equation}
\label{eqn5.8}
u(x,t) = \frac{(x - 14000)\pi}{5400 h(x,t)} \mbox{cos} \left( \pi \left( \frac{4t}{86400} - \frac{1}{2} \right) \right).
\end{equation}

The D2Q9 velocity model is used with $f^{eq}$ defined by Equation \eqref{eqn3.35}. The value of the gravitational acceleration used is between $0$ and $\frac{3}{5}$, i.e $g \in (0,\frac{3}{5})$ and $\lambda = 1$. We choose to use Proposition \ref{prop4} since we have shown in Example 1 that, the two equilibrium distribution functions in Equations (\ref{eqn3.26}) and (\ref{eqn3.35}) behave the same way when modelling shallow water flows. Similarly, two-dimensional code was used to produce the numerical results for a one-dimensional problem. Periodic boundary conditions were used at the upper and lower walls.  

First we need to discuss the stopping criterion and time accuracy of the algorithm. The analytical solution of the flow is known and will be used for validation of the numerical solution. The methodology used in Example 1 will be used in this example. The $L^2$- error norm is used as defined in Equation (\ref{eqn:error}). 
\begin{figure}[ht]
\centering
\includegraphics[width=0.55\textwidth,height=0.3\textheight,angle=0]{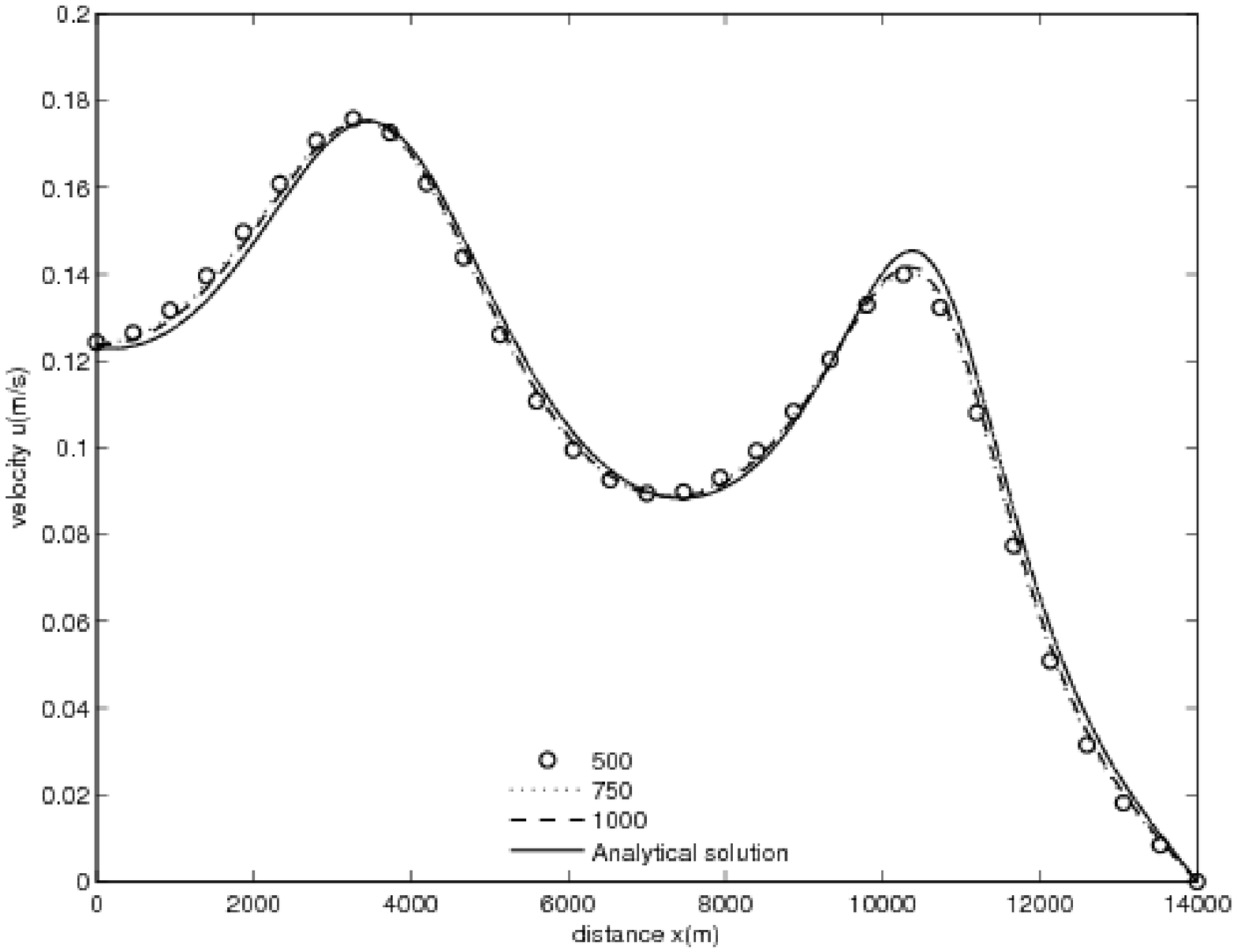}
	\caption{Numerical and analytical free surface for the tidal wave flow at time $t = 9117.5$ s.}
	\label{wave_velocity1}
\end{figure}

\begin{table}[ht]
	\caption{Comparison of numerical and analytical solutions using $L^2$- error norm.}\centering
	\label{error_approximation}
		\begin{tabular}{ c  c }
			\hline\noalign{\smallskip} 
			\textbf{Lattice size (m)} & \textbf{$L^2$- error norm}\\ [3pt]
			\hline \hline																	
			\noalign{\smallskip}			
			\multirow{1}{*}{$\Delta x = 7$}	 	 & $5.27 \times 10^{-2}$\\					
		 	\multirow{1}{*}{$\Delta x = 14$}	 & $6.39 \times 10^{-2}$\\
			\multirow{1}{*}{$\Delta x = 28$}	 & $6.68 \times 10^{-2}$\\									\noalign{\smallskip}													
			\hline\hline
		\end{tabular}
		\end{table}

Three uniform lattices with  $500\times50$, $750\times50$ and $1000\times50$, which correspond to $\Delta x = 28$ m, $\Delta x = 14$ m and $\Delta x = 7$ m, respectively, were used. For the numerical computation, $\tau = 0.6$ and $e = 200$ m/s is used. Similar to Example 1, the value of $\lambda$ was varied between -4 and 7 with $g = \frac{1}{3e}$. It can be observed that the algorithm converged when $t = 9117.5$ s. 

To quantify the results obtained, a comparison of the asymptotic results in Equations (\ref{eqn5.7}) and (\ref{eqn5.8}) are compared with the computed solution. Figure \ref{wave_velocity1} shows a comparison of the numerical solutions with the analytical solution at $t = 9117.5$ s, where $g = 0.0017$ ($\frac{1}{3e}$) was used. It is clear that the results compare favourably. It was found that using $\Delta x = 7$ m gave slightly better results, refer to Table \ref{error_approximation}. It can be concluded that, when lattice size is decreased then accurate results are obtained. Similar behavior has been observed in Example 1. For the water depth it was found that the relative error was about 2 $\%$. The numerical and analytical solutions for the free surface are depicted in Figure \ref{level}. 

\subsection{Example 3: Flow over a sudden-expansion channel}

\begin{figure}[ht]
\centering\includegraphics[width=0.55\textwidth,height=0.3\textheight,angle=0]{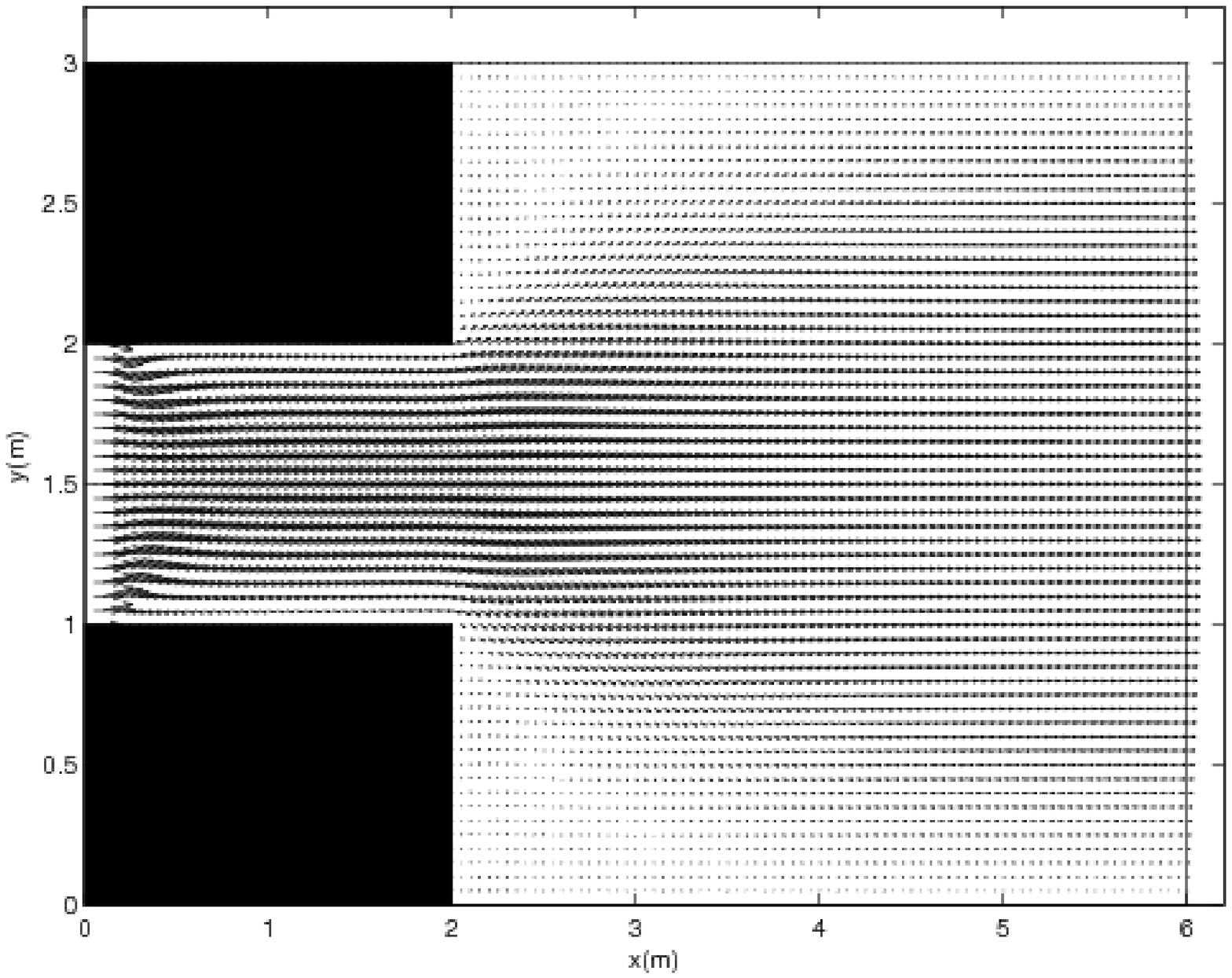}
	\caption{Sudden-expansion channel: velocity field, where $g = 0.15$.}
	\label{exp_vec}
\end{figure}

In this example, we consider a two-dimensional (2D) flow over a channel with a symmetric sudden-expansion. The idea behind this example to simulate circulation flow. The channel expansion ratio is 3:1 with a channel expansion of 3 m wide and 4 m long. The entrance of the channel is 1m wide and 2m long, refer to Figure \ref{exp_vec}. In this case the bed slope and friction at the bottom are neglected.

For numerical computations, the D2Q9 velocity model is used with $f^{eq}$ defined by Equation (\ref{eqn3.35}). The structure of the grid contains $120 \times 60$ lattice points with, $\Delta x = \Delta y = 0.05$ m, $\Delta t = 0.025$ s and $\tau = 1$. As above the speed of the lattice $e$ is given by $e = \Delta x /\Delta t$. The boundary conditions are prescribed as below: the water level $h = 0.16$ m is used at the outflow boundary, zero gradient of depth is specified together with the discharge $q = 0.032$ $m^3/s$ at the inflow boundary. In addition the velocity, $u = 0$, is imposed at the inflow.

Figure \ref{exp_vec}, shows the velocity field with a bit of circulating flows on both sides of the channel. From this figure, it can be concluded that SWEs are capable of simulating circulations that occur in shallow water flows provided the parameters are chosen in order to satisfy the stability notion.

To test the effect of the stability structure the following tests were undertaken. In the first set of simulations, different values of $g$ between $0.001$ and $0.5$ were used in the computation with the parameter $\lambda = 1$ fixed. The steady state solution was reached using the convergence criterion in Equation (\ref{eqn5.9}). The number of iterations required before convergence to steady state was attained are presented in Table \ref{sudden_expansion}. It can be observed that when the value of $g$ was varied in the interval $0 < g < \frac{3}{5e}$, the algorithm converged. 
\begin{table}[ht]
	\caption{The summary of gravity values for $\lambda = 1$ on a $120 \times 60$ grid.}
\centering
	\label{sudden_expansion}
		\begin{tabular}{ c  c }
			\hline\noalign{\smallskip}
			\textbf{Gravity (g)} & \textbf{Number of iterations}\\ [3pt]
			\noalign{\smallskip}\hline \hline
	  	\multirow{1}{*}{0.001}	  & 21645\\																\multirow{1}{*}{0.08}	 		& 13432\\										 \multirow{1}{*}{0.15}	    & 11123\\
			\multirow{1}{*}{0.23}	 		& 31373\\
			\multirow{1}{*}{0.3}	 		& |\\		
			\multirow{1}{*}{0.5}	    & |\\														\noalign{\smallskip}							
			\hline\hline
		\end{tabular}
\end{table}

In the second set of tests, the value of $g$ was fixed and the parameter $\lambda$ was varied between $-2$ and $12$. The iterations required before convergence to steady state are presented in Table \ref{sudden_expansion_lambda}. Similar observations made above can also be noted here.

\begin{table}[ht]
	\caption{The summary of $\lambda$ values using $g = 0.1667$ ($g = \frac{1}{3e}$).}\centering
		\label{sudden_expansion_lambda}
\begin{tabular}{ c  c }
			\hline 
			\noalign{\smallskip}
			\textbf{Values of ($\lambda$)} & \textbf{Number of iterations}\\ [3pt]
			\hline \hline
			\noalign{\smallskip}
	  	\multirow{1}{*}{-2}	  	& 23039\\																					
			\multirow{1}{*}{4}	 		& 23039\\								
		 \multirow{1}{*}{7}	    	& 23039\\
			\multirow{1}{*}{12}	 		& |\\												\noalign{\smallskip}\hline\hline
		\end{tabular}
			\end{table}
			
In conclusion, for this example as well, the parameter values within the ranges prescribed by the stability notion guarantee convergence.
			
\section{Conclusion and Further Work}\label{section5}
A stability structure defined in \cite{Banda} to investigate the stability of the LB equations which are currently being applied to simulate SWEs has been discussed. The models which were chosen were two-dimensional (2D) and have sufficient symmetry, which is a dominant requirement for the recovery of SWEs from lattice Boltzmann equations \cite{Salmon}. In this paper, a stability notion which can be used for constructing lattice Boltzmann equations for SWEs is proved. With the stability requirement, relations of parameters have been derived for parametrized models.

Three examples were used in this work to test the fully discrete LB method. The theoretical results in Section (\ref{Determination of Parameters}) have been tested. The numerical results verify that the stability structure is an appropriate tool for designing the requisite lattice Boltzmann equations for shallow water equations. This applies to both steady state and time-dependent problems.

For further work the consistency of the lattice Boltzmann Equations for shallow water equations will be investigated. In this paper the models generally used in the literature were found to be stable under certain conditions and their consistency has not yet been proven rigorously.  Hence, an important aspect that requires further research is to verify the consistency of the models. This will then complete the numerical analysis for the lattice Boltzmann method for shallow water equations.

\end{document}